\documentclass{article} 

\pdfoutput=1  

\usepackage{fncylab,enumerate,wrapfig,placeins}

\usepackage{amsmath}

\usepackage{amsthm,amsfonts,amssymb, graphicx}  
\usepackage{relsize} 
\usepackage{accents}  

\usepackage{caption}
 \usepackage{textgreek,upgreek,bm}

 \usepackage{titlesec}

\usepackage[usenames, dvipsnames]{color}

\usepackage[pdfa,hidelinks]{hyperref}
\usepackage{geometry}

\usepackage{optidef}
\usepackage{graphics} 
\graphicspath{{figures/}}
 
\usepackage{balance}
 
\usepackage{verbatim}

\usepackage{soul}

 \usepackage{pdfpages}  

\def\urls#1{{\scriptsize\url{#1}}}
 
\newcounter{rmnum}
\newenvironment{romannum}{\begin{list}{{\upshape (\roman{rmnum})}}{\usecounter{rmnum}
\setlength{\leftmargin}{14pt}
\setlength{\rightmargin}{4pt}
\setlength{\itemsep}{1pt}
\setlength{\itemindent}{12pt}
}}{\end{list}}

\def\Ebox#1#2{%
\centerline{\includegraphics[width= #1\hsize]{#2}}%
}

\def\bl#1{{\color{blue}#1}}
\def\rd#1{{\color{red}#1}}
\def\gr#1{{\color{ForestGreen}#1}}

\newlength{\noteWidth}
\long\def\notes#1{\ifinner
             {\tiny #1}
             \else
             \marginpar{\parbox[t]{\noteWidth}{\raggedright\tiny #1}}
             \fi}
             
            \def\notes#1{\typeout{See notes!}}
\def\archive#1{}

\def\rwm#1{\notes{\gr{rwm:  #1}}}

\def\jjm2021#1{\notes{\rd{jjm 2021:  #1}}}

\def\MC{\text{\MC}}

\def\argmin{\mathop{\rm arg\, min}}

\def\argmax{\mathop{\rm arg\, max}}

\newcommand{\field}[1]{\mathbb{#1}}

\def\Re{\field{R}}

\def\eqdef{\mathbin{:=}}

\newtheorem{theorem}{Theorem}[section]

\newtheorem{proposition}[theorem]{Proposition}

\usepackage{cleveref}

\Crefname{corollary}{Corollary}{Corollaries}
\Crefname{eqnarray}{eq.}{eqs.}
\Crefname{equation}{eq.}{eqs.}

\Crefname{figure}{Fig.}{Figs.}
\Crefname{tabular}{Tab.}{Tabs.}
\Crefname{table}{Tab.}{Tabs.}
\Crefname{lemma}{Lemma}{Lemmas}

\Crefname{theorem}{Thm.}{Thms.}
\Crefname{definition}{Definition}{Definitions}
\Crefname{section}{Section}{Sections}
\Crefname{proposition}{Prop.}{Propositions}
\Crefname{assumption}{Assumption}{Assumptions}
\Crefname{example}{Example}{Examples}

\def\bfmd{\bfmath{d}}

\def\bfmg{\bfmath{g}}

\def\bfmg{g}

\def\bfmd{d}

 \def\FRAC#1#2#3{\genfrac{}{}{}{#1}{#2}{#3}}

\def\ddt{{\mathchoice{\FRAC{1}{d}{dt}}%
{\FRAC{1}{d}{dt}}%
{\FRAC{3}{d}{dt}}%
{\FRAC{3}{d}{dt}}}}

\def\clL{{\cal L}}

\def\clT{{\cal T}}
\def\clU{{\cal U}}

\def\cX{c_{\text{\tiny X}}}

\def\cdG{c_{{\text{\lower1pt\hbox{d}}}} }

\def\As#1{\noindent
	\textbf{(#1)}  \ \ }

\usepackage{bbold} 

\RenewEnviron{BaseMiniExclam}[7]{%
	\selectConstraintMult{#1}%
	\renewcommand{\localOptimalVariable}{#2}%
	\begin{subequations}
		\ifthenelse{\equal{#7}{b}}{\allowdisplaybreaks}%
		#4
		\begin{alignat}{5}
		\bodyobj{#2}{#3}{#6}{#5}
		\BODY
		\end{alignat}
	\end{subequations}%
	\setStandardMini
}

\makeatletter
\def\thanks#1{\protected@xdef\@thanks{\@thanks
		\protect\footnotetext{#1}}}
\makeatother
\title{\Large \bf
	Reliable Power Grid: Long Overdue Alternatives to Surge Pricing
}

\author{Hala Ballouz$^{1}$,
Joel Mathias$^{2}$, 
Sean Meyn$^{2}$, 
Robert Moye$^{2}$, 
and Joseph Warrington$^{3}$
\thanks{$^*$Funding from the National Science Foundation under award EPCN 1935389.   Thanks also  to funding from  the State of Florida,  through a REET (Renewable Energy and Energy Efficiency Technologies) grant. Many thanks to Prof.\ Frank Kelly who suggested we survey the history in telecommunications economics to investigate parallels with the power industry. The Simons Institute, Berkeley, is gratefully acknowledged for hosting SM and JW in Spring 2018
 }
\thanks{$^{1}$
President, Electric Power Engineers,
Austin, Texas}%
\thanks{$^{2}$Department of Electrical and Computer Engineering, University of Florida, Gainesville, FL 32611}%
\thanks{$^{3}$Home Experience LLC (UK), Cambridge, UK}%
}

\date{March 12, 2021}
 
\begin{document}
 
 \includepdf[page=1]{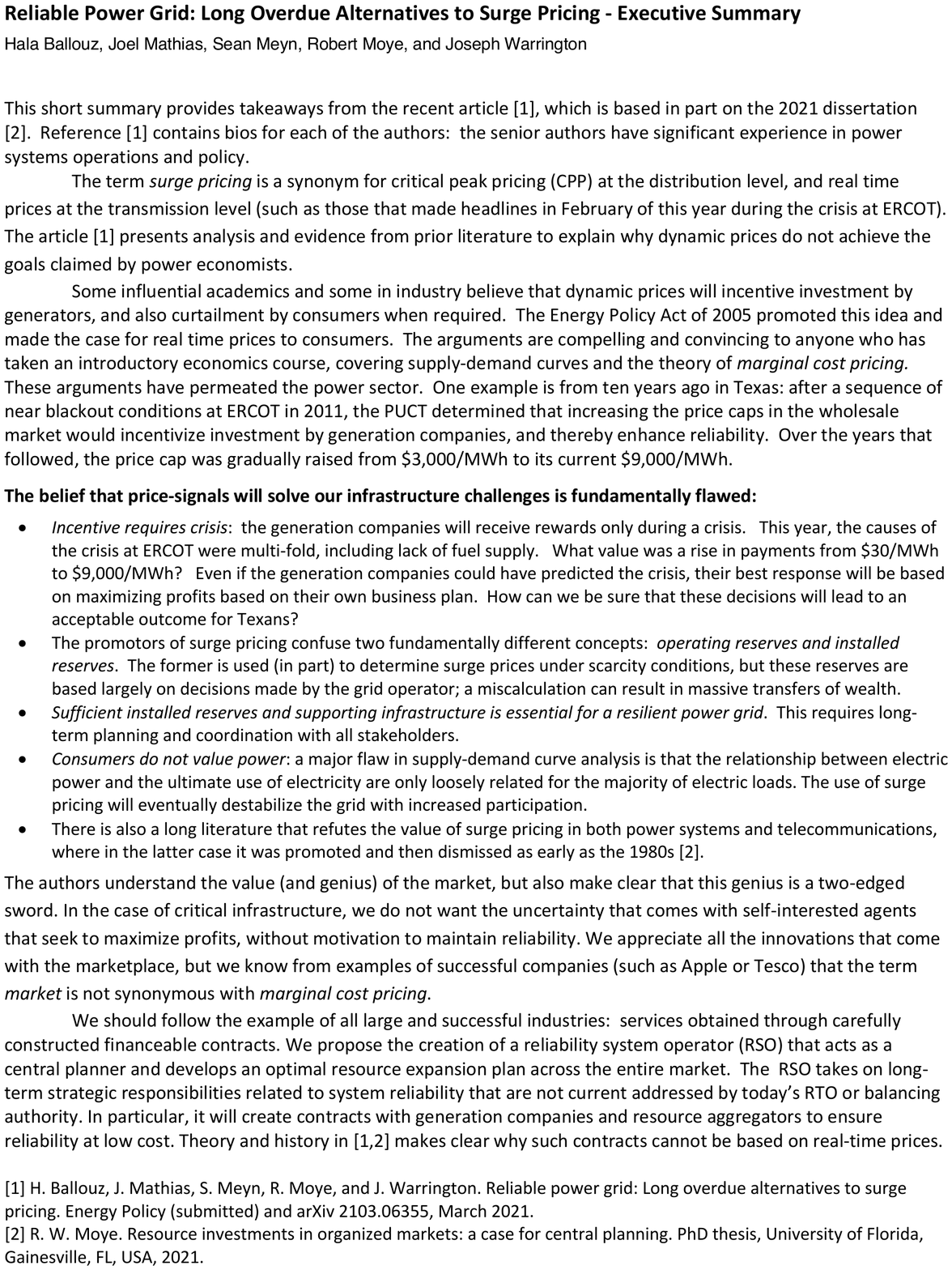}

\clearpage

\maketitle
\thispagestyle{empty}

\begin{abstract}
 
Since the 1970's, it has been recognized that demand-side flexibility of electric loads can help to maintain supply-demand balance in the power grid.  One goal of the  Energy Policy Act of 2005 was to accelerate the adoption of demand-side resources through market design, with competitive equilibrium theory as a guiding philosophy.  This paper takes a fresh look at this economic  theory that is the motivation for pricing models, such as critical peak pricing (CPP), or surge pricing, and the demand response models advocated  in the Energy Policy Act.      

The economic analysis in this paper begins with two premises:  1) a meaningful analysis requires a realistic model of stakeholder/consumer rationality, and 2) the relationship between electric power and the ultimate use of electricity are only loosely related in many cases.  
The most obvious examples are refrigerators and hot water heaters that consume  power intermittently to maintain their temperature within predefined bounds.   A dynamic economic model is introduced based on these premises.  This is used to demonstrate that with the use of CPP and related pricing schemes will eventually destabilize the grid with  increased participation.
  Moreover,  it is argued that the optimal dynamic prices (supporting a competitive equilibrium for the dynamic model) do not lead to a robust control solution that is acceptable to either grid operators or consumers.  
These findings are presented to alert policy makers of the risk of implementing real-time prices to control our energy grid.  
 
Competitive equilibrium theory requires a coarse description of a real-world market, since complexities such as capital costs are not included.
The paper explains why these approximations are especially significant in the power industry.  It concludes with policy recommendations to bolster the reliability of the power grid, with a focus on planning across different timescales and alternate approaches to  leveraging demand-side flexibility in the grid.


\end{abstract}



 
%
%
  \clearpage
 
\section{Introduction}
\label{s:intro}

Several market models in the United States, particularly in the Electric Reliability Council of Texas (ERCOT), were structured around energy settlement through real-time pricing as a mechanism to drive reliability. These models are mostly prevalent as energy settlement markets on the bulk power systems.   However, there is evidence that they do not offer the foundation to drive a reliable power grid.  The trends over the past decade suggest that similar market structures will be imposed on the distribution system in order to leverage demand-side flexibility as a commodity. This paper explains why policy makers should find this trend alarming and, in particular, sets out to explain troubling faults in the theoretical underpinnings of marginal cost pricing.

 
 The term  \textit{surge pricing}  is a synonym for critical peak pricing (CPP) at the distribution level, and real time prices at the transmission level (such as those that made headlines in February of this year during the crisis at ERCOT).   


\subsection{Why not surge pricing?}  

In the Spring of 2018, the Simons Institute for the Theory of Computing hosted a program on \textit{real-time decision making}.  The April workshop on \textit{New Directions in Societal Networks},  which included the head of market operations at Uber and the director of product management at Lyft as speakers, was one of the high points of the program.   Both these speakers discussed the challenges in the design of pricing mechanisms, and both discussed the \textit{failure} of surge pricing.    There are many reasons for failure, but one in particular should resonate with policy makers in the energy industry:  \textit{surge pricing does not achieve the goal predicted by theory},
 due to significant gaps in the underlying assumptions.
Among the explanations is  delay:  if drivers receive news of a significant high price event,   they will likely predict that prices will have stabilized by the time   they reach the congested region.    

It is well known in control theory and every other domain of ``real-time decision making'' that a carefully designed feedback system may not behave as expected in the presence of delay.   This principle is the starting point of \cite{roodahmit10a,roodahmit12,seel2018impacts}, concerning the potential instability of the power grid with the introduction of real-time prices. 

With this background in mind, we consider  price dynamics in the power grid.

\jjm2021{could be moved as a lead-in to contributions or deleted: Surge pricing and markets driven by the marginal cost of power, particularly in deregulated markets with siloed stakeholders, fail to meet reliability standards pertaining to the delivery of energy to the consumer.}

\paragraph{Supply side dynamics}
The extraordinary events in the ERCOT region during February 2011 and again in February 2021 provide a concrete example of how surge pricing may not provide the intended incentives; instead, such pricing undermines the reliability of the grid, while simultaneously causing tremendous financial strain to the stakeholders.    
In 2011,  wholesale electricity prices in Texas hit the price cap of \$3,000/MWh for several hours in February~\cite{HOUSTON2011},   and recurring high price events occurred on a daily basis for several weeks during the summer.      Remedies were proposed following an investigation by the Public Utility Commission of Texas (PUCT),  and one of these was to raise the price cap, driven by academic research  centered around the marginal value of power \cite{hogan2013electricity,surendran2016scarcity}.  The following is taken from the introduction of the  2012 Brattle Group report \cite{newell2012ercot}:   (the PUCT)  \textit{has implemented
a number of actions to ensure stronger price signals to add generation when market conditions
become tight. The PUCT has enabled prices to reach the current \$3,000/MWh offer cap under a
broader set of scarcity conditions and is considering raising offer caps to as high as  \$9,000/MWh, among other measures.}


The key conclusion in   \cite{hogan2013electricity} is that ``\textit{suppressed prices in real-time markets provide inadequate incentives for both generation investment and active participation by demand bidding}''.
This conclusion, at face value, would be correct if not for unavoidable delays, both short-term and long-term (amongst other `imperfections'):  it takes time to create capacity (years for most resources and transmission lines),  to make decisions to switch circuits under emergency conditions, or even to start up a thermal power plant.\jjm2021{to SPM on Mar 10: this is said again below, so removing: it also requires much more certainty into the prediction of those future surge prices to justify long term expensive assets to be built.} This limits the effectiveness of surge pricing, exactly as in the ride sharing business.  In certain markets like ERCOT, the uncertainty introduced by such mechanisms diminishes incentive to plan and build new power plants for the future, or properly maintain assets. 
More than 30 Texas power plants owned by Fortune 500 companies failed in the 2011 winter freeze and again in 2021, ``despite warnings about the need to winterize their equipment''~\cite{WP_ERCOT_Mar6_2021}.

Short-term delay played a role in the   2021 power systems crisis at ERCOT:   ``\textit{Some of these units on outage were likely unable to secure gas on such short notice as the weekend’s gas prices rapidly increased from \$7 per mmBTU on  Thursday to \$150 per mmBTU by the weekend due to supply concerns around freeze-offs and heating demand.}''\footnote{\urls{https://www.woodmac.com/news/editorial/breaking-down-the-texas-winter-blackouts/full-report/}}


\paragraph{Demand side dynamics}
The issue of delay is far more exotic when we consider demand-side participation, especially in the case of residential consumers. An examination of user preferences is required to make this precise. 
 
In our economic model, we assume the following features of a so-called ``rational agent'' at his or her home:
\begin{romannum}
\item  The refrigerator temperature should remain within prescribed bounds;
\item  A hot shower should be hot, but not too hot  (the thermostat setting should be respected);
\item A pool should be cleaned regularly (say, 4-10 hours per day and  40 hours per week);
\item  Subject to the above ``quality of service''  (QoS) constraints,  the  power bill should be minimal.
\end{romannum}
The acronym QoS is borrowed from the telecommunications literature, and is used here to emphasize a common concern in these disparate applications.  

\begin{figure}[h]
	\Ebox{.85}{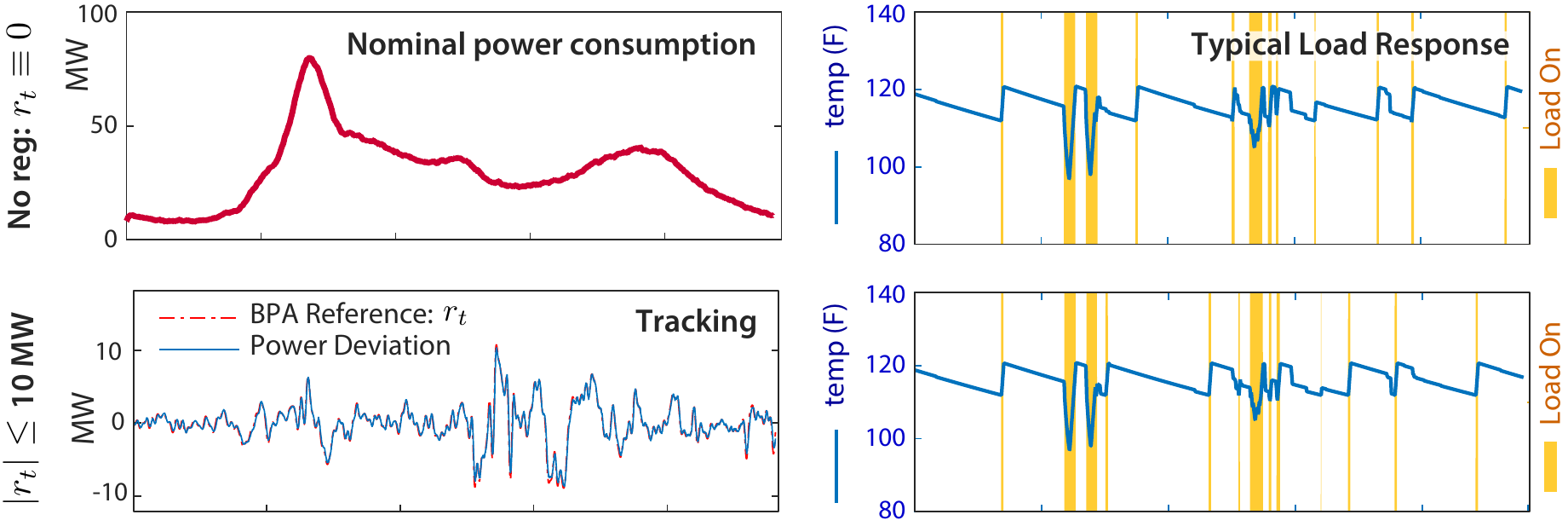} 
	\caption{The first row shows nominal power consumption from a collection of 100,000 residential water heaters,   and  the behavior of a single water heater.    The second row shows how these water heaters can track the Bonneville Power Authority (BPA) balancing reserve signal.    The behavior of the individual water heater shows no apparent changes. }  	
 	\label{f:WHTemp1}
\end{figure}

It may surprise many readers that a standard residential water heater consumes power quite rarely:
during a typical day, power consumption is roughly periodic, with periods typically ranging from two hours to ten hours. As an example, a water heater with a six hour period will consume power for only five to ten minutes in that period.  The first row of \Cref{f:WHTemp1} shows behavior of a typical residential water heater under heavy use (explanation of the plots on the second row is postponed to \Cref{s:sanity}). 
Few residents are aware of these power consumption patterns,  which is a great benefit.   As long as hot water stays within bounds,  
the consumer is content with the outcome. It is not important to the consumer precisely when power is consumed.

Similar statements can be said about many loads:  a residential refrigerator has a shorter duty cycle,  and more symmetric power consumption cycle. Electrical demand for pool cleaning or irrigation is often substantial, and yet also highly flexible.

\begin{wrapfigure}[13]{r}{0.375\textwidth}
 \centering 
    \includegraphics[width=0.95\hsize]{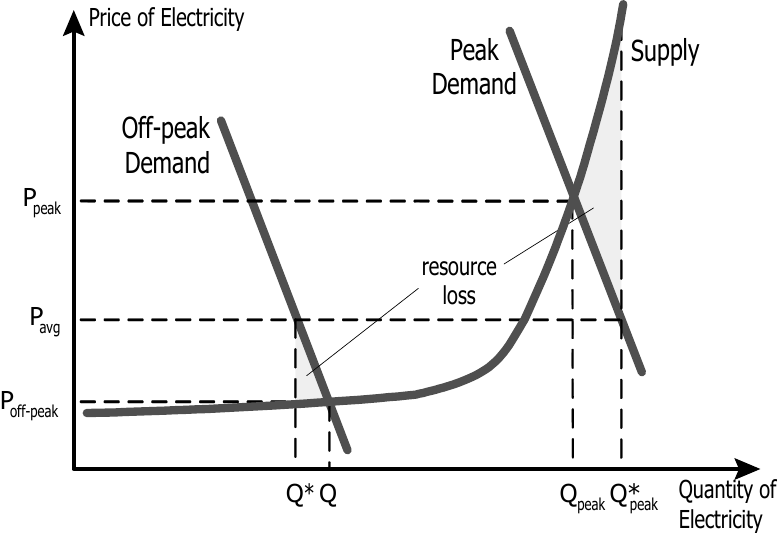}
	\caption{Is demand for electricity   a continuous function of price?    }  	
 	\label{f:DOE_S=D}
 
\end{wrapfigure}

Many government reports and academic articles currently adopt  an entirely different model of user preferences.    It is typically assumed that   power consumption is a continuous function of price,  and consequently price signals can be used to control the grid. 
 The plot shown in \Cref{f:DOE_S=D} is adapted from Fig.~B.1 of the DoE report
 \cite{qdr2006benefits} is an illustration of this assumption.   Similar plots appear in the aforementioned paper \cite{hogan2013electricity}, and  as Fig.~8.3 in \cite{borenstein2005time}.
 The analysis in  \cite{borenstein2005time} takes for granted that power demand changes smoothly with price, and that the correct price is the marginal cost of generation;  conclusions that are inherently flawed in many cases, as we explain in this paper. 
The assumption that power consumption varies smoothly is a starting point in many other policy papers and academic articles  \cite{chomey10, wankownegshamey10,kizman10b,hogan2013electricity,zavani11a,jor19}.  Based on this assumption, it follows that the grid can be successfully managed with  dynamic prices to consumers.

A goal of this paper is to make clear why price response cannot be predicted through the standard supply-demand curve shown in 
\Cref{f:DOE_S=D} (in most cases),  and the reason will sound familiar since it is rooted in temporal dynamics (much like the delayed response of a driver for Uber or Lyft).

\begin{figure}[ht]
	\Ebox{.65}{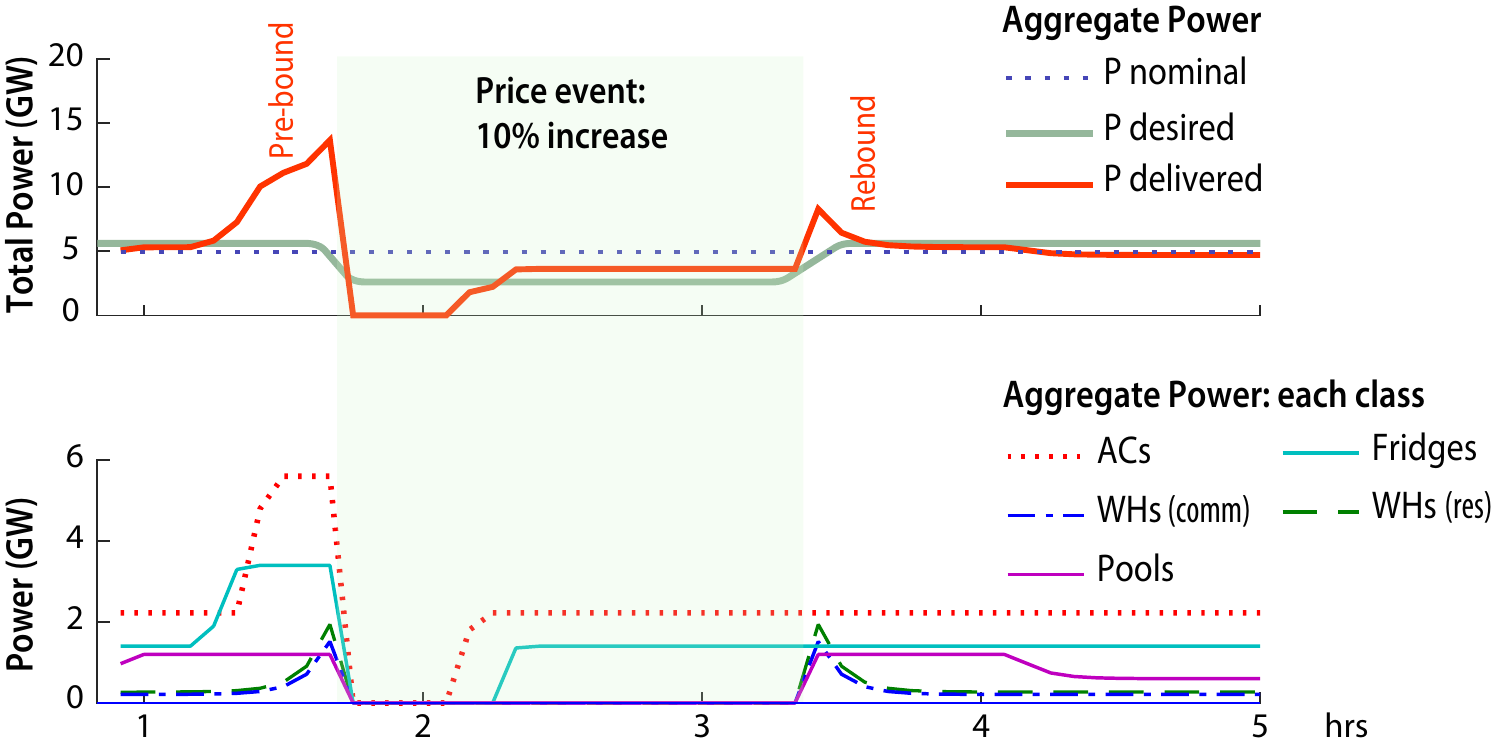}
	
	\caption{Optimal response to a 10\%\ increase in price over ninety minutes \cite{cammatkiebusmey18}. }
	\label{f:boom} 
\end{figure}

\jjm2021{not needed or meaningful: In an alternative model to surge pricing, }Thermostatically  controlled loads (TCLs)   such as water heaters and refrigerators can be regarded as energy storage devices.     It is shown here that this characteristic leads to price responses that are highly dynamic and in general discontinuous as a function of price.    \Cref{f:boom} shows an example of how residential and some commercial loads would respond optimally to a CPP event, which includes a 13~GW (gigawatt) drop in net-load when all participating loads turn off simultaneously at the onset of the price event.
Full details of this experiment are contained in    \Cref{s:idiots}.

   We also demonstrate how we can take advantage of realistic models of user preferences to create valuable resources to balance the grid.  In the case of an electric water heater or refrigerator, the ``user preference'' is found on the thermostat settings.

\subsection{Some History and Recent Developments}
\label{s:history}



There is a history that makes clear why we should seek alternatives to   price signals   in order to better control the grid. 

The foundations of price-based demand response are rooted in the theory developed by Dupuit in the 1840s, Hotelling in the 1930s, Vickrey in the 1950s, 
and Schweppe in the 1980s \cite{dupuit1844, hotelling1937, vickrey1955, schweppe1988}. All of these researchers looked at the relationship between the value to a consumer and the marginal cost to provide the service. Whether this was related to the use of a water system, railroads, bridges or electricity, the focus was on the ``demand'' for the service or product. However, in contrast to the direct value to the consumer for moving a product on a railroad, or the value of crossing a bridge versus taking longer alternate routes,  electricity has always presented a challenge because the product that the customer values is not electricity per se, but the heating, cooling, hot water, or lighting that comes from it. This disconnect fundamentally changes the prices-to-loads compact that has been the subject of much research over the past few decades. 

We have witnessed the \textit{prices to devices} debate before in the telecommunications industry.  Throughout the 1980s, it was argued that real-time prices could help manage telecommunication traffic.     Following deregulation of their industry in the 1980s, telephone companies investigated real-time pricing (what they called ``measured rates'') for local service. At the start of this effort, the assumption by economists and policy makers was that measured-rate pricing could ``substantially increase economic efficiency'' \cite[pg 2]{park1987}.
However, in both study and practice, and looking specifically at the net welfare effects on their customers, flaws were found in this way of thinking: ``[C]ontrary to conventional wisdom, measured rates will at best produce very modest efficiency gains; more likely, they will result in small efficiency losses \cite[pg 3]{park1987}.''   In a statement that was surprising at the time,   the authors concluded that ``if local measured service is desirable public policy, it must be justified on grounds other than economic efficiency.''

In the mid-1980s, several state regulatory commissions had approved various measured-rate pricing plans for telephone companies, and many more proposals were pending. Further, many argued that this was the only way to manage the Internet.  Studies revealed that the real-time tariffs  would be 1) ``too complicated for telephone subscribers to cope with \cite[pg 6]{park1987},'' and 2) the systems required to implement them would also be too complex. 
The pricing debate closed before the end of the decade, with the successful implementation of distributed and automatic control systems that are now also found in every computer and cell phone  \cite{lobluhinmey19}.

With the integration of renewable energy resources, which are characterized by high fixed costs and  almost zero variable costs (characteristics similar to the telecommunications industry), the research and conclusions reached by the telecommunication industry are worth consideration by policy makers in the power sector.


\smallskip

These issues are particularly prescient in demand response programs, which fall broadly into two categories: (i) incentive-based programs with a focus on engineering design, e.g. \textit{direct load control} (DLC) or \textit{interruptible load programs}  (IBP), and (ii) price-based programs with a focus on market efficiency, e.g., \textit{real-time pricing} (RTP), \textit{critical peak pricing} (CPP), and \textit{time-of-use} (TOU) tariffs --- see \cite{albadi2008summary} for a detailed survey.   

There is a substantial literature exploring the potential difficulty with real-time prices in the power sector.  
Empirical studies, such as those surveyed in  \cite{qdr2006benefits}, show that a period of  high prices often induces a response from consumers, but it is often not the desired response.  Most commonly noted  is the ``rebound'' effect, in which power consumption increases dramatically following a period of time during which prices are high  \cite{qdr2006benefits,Luetolf2018ReboundEO}.


Moreover, studies have found that price-based demand response has a disproportionately negative impact on indigent, disabled, and elderly consumers, as {CPP} and {TOU rates} can be a punitive signal to these consumers to reduce their power usage, even at a detriment to their quality of life \cite{ale10}. The debacle of the Houston-based company Griddy, 
which is a retail power provider that charges  customers based on wholesale electricity prices, is a cautionary tale: following the winter storm in February 2021, many customers were left with bills amounting to thousands of dollars, leading to multiple lawsuits as well as ERCOT revoking its access to the Texas grid \cite{Griddy}; one of many detrimental financial impacts of surge pricing during that event.

Research from social and economic theory suggests that price-based incentives tend to discourage civic responsibility \cite{titmuss1970gift, frey1997cost}. Formulating demand-response programs  based solely on economic incentives, especially real-time prices, can be deeply problematic \cite{sovacool2014diversity, he2013engage}: consumers tend to value comfort, autonomy, ease-of-use, and privacy over financial benefits \cite{xu2018promoting, parrish2020systematic}.  A study investigating the attitudes of residential consumers in Great Britain towards participation found that automated direct load control (with an override feature) was more popular than TOU or dynamic pricing, because consumers perceive the former as less complex and providing greater flexibility and autonomy as compared to the latter \cite{fell2015public}.

Finally, we recommend the review articles \cite{spe16naive,spe18naive} for more history on re-regulation of the power industry,  and a different perspective on ``\textit{...how
and why energy markets can never resemble the idealized markets}''

\subsection{Contributions} 

The technical contribution of this paper is the introduction of a dynamic competitive equilibrium model, based on a realistic model of consumer preferences.    The model illuminates a flaw in common assumptions regarding real-time prices:   the  equilibrium price has little to do with marginal cost or marginal value,   and standard pricing models can drastically reduce system reliability.   

These technical results alone demand a rethinking of marginal cost pricing.   The authors' experience in the power sector  (three decades each for two of the authors)   inform a broader set of  policy recommendations.   Below is a summary:

\paragraph{Planning}
  A reliable grid requires planning to address  challenges across different time-scales: building capacity for reliability    
  requires long-term planning spanning years; resource allocation using load forecasts occurs over time-scales of several days; whereas rejecting disturbances caused by uncertainty of wind or solar generation, or for that matter, disruptive consumer  behavior, needs to be addressed on time-scales of a few minutes. Planning under so much complexity and uncertainty requires cooperation among experts in concert with all stakeholders. We cannot hope to achieve our long-term reliability goals through short-term marginal cost pricing models.

  As in most other industries, services should be obtained through carefully constructed financeable contracts.  We propose the creation of a \textit{reliability system operator} (RSO)
  that acts  as a central planner and develops an optimal resource expansion plan across the entire market.
The RSO takes on many of the responsibilities of today's   RTO or balancing authority (BA).  In addition, it will create contracts with generation companies and resource aggregators to ensure reliability at low cost. The theory in this paper explains why such contracts cannot be based on real-time prices.
More on the responsibilities of the RSO is contained in  \Cref{s:RSO}.


An example of effective planning and use of contracts can be found at PJM.  They understood that the most severe winter generation capacity shortages may be overcome using regulated forward supply contracts, along with stiff penalties for failure to deliver the contracted electric capacity at times of tight supply.     PJM's three-year-ahead system requires utilities to procure capacity to cover their customers' aggregate demand\footnote{\url{https://learn.pjm.com/three-priorities/buying-and-selling-energy/capacity-markets.aspx}}.   Another example is the British capacity market,\footnote{\url{https://www.emrdeliverybody.com/cm/home.aspx}} which runs annual auctions for both single-year payments a year ahead of time, and 15-year payments four years ahead of time. The purposes of these auctions are respectively to ensure winter adequacy in the short term, and to de-risk the building of new generation in the long term.  
Coordinated mechanisms such as these deploy the necessary time-scales and service contracts that honor the ``delays'' that are otherwise ignored in the current simplified textbook assumption that rational generation companies already have the required incentives to hedge financially against any future inability to deliver power.


%
%
%

\paragraph{Distributed intelligence}
  
	To obtain reliable grid services from flexible loads, we can follow the success story of the Internet, where supply-demand balance is achieved through distributed control across the network.   It is possible to design decision rules so that strict performance guarantees are met, for both the grid and consumers,  without sacrificing privacy.
	
	Through a robust, decentralized, automated control design, the emerging science of \textit{demand dispatch} can serve the needs of the grid operator while ensuring consumer-side comfort, privacy, and ease-of-use \cite{brolureispiwei10,chehasmatbusmey18, matbusmey17, matkadbusmey16}.   Resource aggregators  (e.g., Enbala Power Networks, Comverge, CPower, Enel X, etc.,   and perhaps a broader role for utilities) 
are required to implement any demand dispatch design.  There is a need for economic models and regulatory frameworks that address the role of these aggregators,   with grid reliability among the highest priorities   \cite{lu2020fundamentals}.

Policy makers should work with engineers to strengthen mandatory standards for appliances,\footnote{One such example is CEA 2045-2013 (ANSI): Modular Communications Interface For Energy Management} so that they can provide \textit{virtual energy storage} for   reliable grid services, without negatively impacting the consumer.

%

%
%


\subsection{Organization}

The technical material in this paper is intended to demonstrate how marginal cost pricing does not reflect the reality of reliable power production, transmission, distribution, and consumption. This requires background on both grid operations and microeconomics.

 \Cref{s:models} is intended to ``set the stage'' for a dynamic economic model, focusing on the needs of   three ``players'':   the  balancing authority (BA),   the generation companies, and the consumers of electricity.   
We recall the need for balancing resources,   and how these might be obtained by exploiting the  large gap between the needs of the consumer and the needs of the BA.      The flexibility of a load is highly dependent on its intended use,   and the \textit{time-scale} of flexibility.


\Cref{s:idiots} reviews dynamic competitive equilibrium theory, and explains why the price signals proposed today will eventually destabilize the grid. 
This conclusion holds even under ideal conditions---no communication delay or other imperfections.   A price surge will cause instability if consumers are permitted to  optimize based on their preferences (such as hot water at a low price).    The optimal outcome is also described, which corresponds to the efficient outcome in a competitive equilibrium model.  
Conclusions, policy recommendations,  and suggestions for future research are contained in \Cref{s:alt,s:con}.

\section{Agent Models}
\label[section]{s:models}

The electric grid is an interconnected network of transmission (high voltage) and distribution (low voltage) lines and other devices that deliver electricity from generating resources to electric load consumers across the network.  Three agents are involved in the transactional process of balancing energy delivery through this network, and are used in the competitive equilibrium model introduced in this paper.

\subsection{The Balancing Authority}
\label{s:grid}

Management of electric power in most industrialized countries is based on a decomposition of the power grid into geographic entities known as balancing areas. The BA is responsible for grid reliability in a balancing area. Each BA manages resources to balance power supply (generators) and demand (consumers), and regulate the flow of power among neighboring BAs.

\begin{figure}[h]
\Ebox{.85}{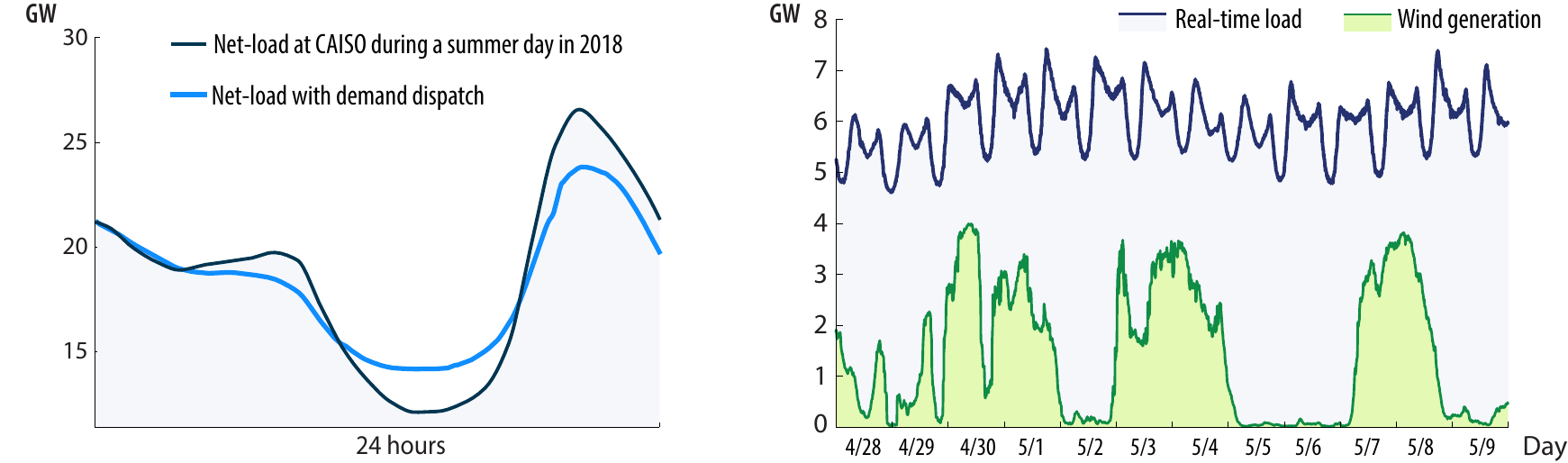}

\caption{Energy from the sun and wind bring volatility.  The plot on the left shows a net-load at CAISO during a very sunny day in 2018.  Shown on the right is a typical week at BPA. }
\label{f:WindAndSun} 
\end{figure}

\Cref{f:WindAndSun} illustrates the tremendous balancing act faced by balancing authorities in the U.S.\ today.   Shown on the left is the net-load  (load minus  generation from renewables) during one day at CAISO in 2018.  Observe that the minimum occurs early in the afternoon when the sun is at its peak,   and then the peak in net-load occurs soon after dusk.  The ramp from minimum to maximum is over 15~GWs,   which is unlike anything seen before 2010.   The load and wind generation shown on the right hand side are taken from a typical week in the Bonneville Power Authority (BPA) in the northwestern U.S.   The ramps in power from wind generation are   significant.

The responsibilities and needs of the BA and utilities in supporting the electric grid and balancing the energy are broad:
\begin{romannum}
\item  Peak shaving and valley filling;
\item Resources that can ramp up (or down) quickly due to a surge in demand or supply of electricity;
\item Black start support;
\item   Addressing high frequency volatility, managed today by selected generators as part of \textit{automatic generation control} (AGC);
\item     
 	Voltage control is achieved using generators at the transmission level;  transformers and other devices manage voltage at the edges of the network;
	
\item    Monitoring the network elements and ensuring that the transmission and distribution lines and supporting devices are adequate to allow the successful delivery of generation to consumers.

\end{romannum}

The BA sends signals to each generator on the grid to ramp up or down its power output in order to balance (or follow) the net load.  Energy storage, particularly batteries,  are being added to the mix of resources available to BAs for balancing energy, with widespread utilization intended to smooth out the intermittencies introduced by renewable resources. These \textit{ancillary services} (regulation, ramping and contingency reserves) are growing in importance as the power grid grows in complexity.

 {Virtual energy storage} (VES), which refers to the aggregation of load-side resources that can provide energy flexibility, will be an important resource for providing ancillary services in the future. For example, Florida Power and Light (FPL) engages 750,000 households in Florida to allow for direct control of water heaters and pool pumps to shed load in case of emergencies. Austin Energy’s Power Partner program allows for brief adjustment of the temperature of the enrolled NEST thermostats by a couple degrees in the afternoon during the hottest days of the year (users can manually override the settings). A larger and older program in France allows the utility to power on water heaters at night. This is beneficial since it is not feasible to turn off nuclear power plants, and the energy they produce has to be consumed.
\textit{Virtual power plants} are defined analogously:  this could be in the form of an aggregation of generators or batteries at the end consumer level,   designed to serve as balancing resources to BAs.

The electric grid, consisting of transmission and distribution wires (and supporting devices),  is the transportation vehicle by which energy delivery is achieved and balanced, and adds another significant layer of complexity.    All of this demonstrates that  maintaining system reliability is a complex task that is not just about the balance of energy.

\subsection{Generation Companies}
\label{s:gens}
 




Generators play a significant role in the reliability of power systems, not only in providing a source of energy to meet the consumer loads, but also to  maintain the stability of the network.  Power systems have to respond to generation  and transmission contingencies as well as load fluctuations in real time.\footnote{During 2019, PJM paid generators roughly \$100 million for regulation services, and \$32 million for synchronous reserves [PJM SOM Report for 2019].}
  Natural Gas and coal are the most prevalent generators in use in the U.S., although renewable generation is claiming a larger share every year. 

Generating resources  that provide grid services beyond energy are compensated through the ancillary service markets. 
Ancillary service products vary  among the organized markets because they have different load and resource characteristics (e.g., one market may have more renewable resources than another), and in general, they may have different reliability needs. For example, CAISO uses a 5-min ramping product, MISO a 10-min ramping product, and ISO-NE does not currently operate a market for a ramping product.  With the increasing penetration of renewables, new system features arise and may create a need for new ancillary services.    

Because energy and ancillary services share the capacity of generation/transmission resources, they must be co-optimized to achieve their most efficient allocation,
 leading to the security constrained economic dispatch (SCED) optimization problem that simultaneously determines the energy and ancillary services designations for resources.
 SCED must also take into account the availability of the network transmission and distribution lines and devices to allow reliable delivery of generation to load.

  
A realistic market analysis must consider all costs associated with generation, and these costs go far beyond  incremental or \textit{marginal} fuel costs:
\begin{itemize}

\item Dispatch costs:   variable  Operation and Maintenance (O\&M) and variable fuel transportation costs;
\item Commitment:  start-up (fuel + maintenance) and no-load energy (fuel);
\item Availability:  fixed O\&M and fixed fuel transportation costs;
	 
\item  Capital costs:  initial investment, renewal and replacement;
\item \textit{Externalities},  including environmental, technological, political, and other factors that are not always directly quantifiable.
\end{itemize}



\jjm2021{comment by Bob, I think: [Sean, I have to dig up the below costs, but wanted to make sure you were going to use this before I do that work.]\\As an example of the costs that the markets incurrs for ancillary serivecs, the PJM market in 2020 alone paid generators \$XXX for regulation, \$YYY for spiniing reserves, and \$ZZZ for non-spinning reserves. In total, this market spent over \$KKK for ancillary services.}

\subsection{Consumers}
\label{s:consumers}

We return here to the consumer, and how a typical household can help to provide VES that can provide many essential grid services.  
The reader can consider their own inventory of devices and preferences.   Here is an incomplete list of residential loads,   organized according to their value for VES capacity as well as the potential risk faced by the consumer.     The term \textit{quality of service} (QoS) refers to the alignment of a consumer's desires with what is offered by a particular device.   
\begin{romannum}
\item  Electric lights,  television, radio and desktop computers.   For these electric loads, the QoS vanishes when the power goes out, so they have no value for grid services outside of a true emergency requiring controlled blackouts.  

\item  Electric vehicles (EVs).    The value of flexibility is tremendous since charging consumes a great deal of power,   and at first sight, there appears to be great flexibility.  However, consumer risk is significant for several reasons:  what is the cost of the additional charge and discharge from the VES?   How to be sure a vehicle will be available when needed in an emergency?     

\item  Thermostatically controlled loads (TCLs),  including refrigeration, water heating and  HVAC (both residential and commercial).    These are energy storage devices that have great potential for VES.   QoS is essentially defined by the thermostat,  which is a user input or set by the factory.

\item  Pool cleaning is a significant source of energy consumption in California, Texas and Florida.    This is obviously highly flexible, as recognized by Florida utilities who contract with homeowners to shut down pool pumps in emergencies.    Commercial irrigation and municipal water pumping require significant energy and are
similarly flexible. 
\end{romannum}
\jjm2021{Mar 10 the latter half of the statement is not clear, not sure anything is added by the first part:\textit{Mathematical modeling of demand dispatch} is addressed here,   as this paper focuses on proving principles in pricing as they propagate to the distribution system and on to the consumer.}

Part of the science required to enable demand dispatch includes models for an aggregation of loads, designed for the dual purpose of modeling QoS and ``virtual state of charge,''  since we interpret the aggregate as a ``virtual battery''.   This language is meant to stress our goal to create VES.    A \textit{virtual battery model} for a collection of thermostatically controlled loads (TCLs) was introduced in  \cite{haosanpoovin15}, and a similar model is proposed in    \cite{chebusmey14,meybarbusyueehr15} for residential pools.   

 The state of charge (SoC) $x_i(t)$ of the $i$th \textit{load class} (i.e., an aggregation of similar loads) at time $t$ amounts to a QoS metric for the population, assumed to evolve according to the linear system,
\begin{equation}
\ddt x_i(t)   = -\alpha_i x_i(t)   +d_i(t),    \qquad 1\le i\le M\,,
\label{e:SoC_ODE}
\end{equation}
in which $d_i(t)$ is power deviation at time $t$,  and $\alpha_i$ is a ``leakage'' parameter.  
For TCLs, the SoC is the normalized thermal energy stored in the population.   In the case of 
swimming pools, it is  discounted-average deviation in desired pool cleaning hours.

The dynamical equation \eqref{e:SoC_ODE} for a collection of TCLs is obtained by simply averaging the temperature of the population, and using standard physics-based models for each individual.   In the case of pool cleaning this emerges as a surrogate for weekly constraints on  hours of cleaning.   Each of these loads is discrete, in the sense that power consumption takes on one of only a few values.     With many hundreds of thousands of loads in play, it is not difficult to approximate the ensemble with a continuously variable input signal $d_i(t)$ for each load class $i$.

The state variable in \eqref{e:SoC_ODE} models only the average QoS deviation,  and hence bounds on $x_i(t)$ represent only a \textit{necessary} condition that QoS constraints are respected for all loads of class $i$.    There is now a mature collection of   distributed control techniques available to guarantee that the necessary condition is also sufficient---a short survey is contained in \Cref{s:sanity}.

\section{Reliability Risks Associated with Dynamic Prices}
\label{s:idiots}

\Cref{f:DOE_S=D} illustrates a component of competitive equilibrium (CE) theory of micro-economics \cite{maswhi1995}, which is briefly surveyed below.
  A dynamic CE model is proposed, from which we obtain our main theoretical conclusions:  (i)  When we adopt a realistic model of consumer preferences, CE prices look nothing like what is predicted by power economists, and (ii)  CPP combined with increasing automation will result in significant reliability risks.
 
 \subsection{Competitive Equilibrium Models}
 
  In the standard CE model, it is assumed that consumers have a utility function $\clU_D$ for electric power (or, in  a discrete time model,  power is substituted for energy consumed in a given time interval).     Subject to this model of rationality, a   consumer at time $t$ will then choose power $P_t $ as the solution to the optimization problem
\begin{equation}
P_t = \argmax_P \{ \clU_D(P) - \varrho_t  P \}
\label{e:CE}
\end{equation} 
where $\varrho_t$ is the price at time $t$. 
  From basic calculus we find that the price is the marginal value of power:
\[
\varrho_t =  \tfrac{d}{dP} \clU_D \,  (P_t)  
\] 
This model of an individual or aggregate of consumers is the starting point of analysis in  \cite{borenstein2005time,chomey10, wankownegshamey10,kizman10b,zavani11a,hogan2013electricity}.
  
A similar calculation shows that the price must   coincide with \textit{marginal cost}  of generation, and from this we obtain the classical supply=demand formula of micro-economics.   It is recognized that the meaning of  marginal cost in power systems is unclear  (see lengthy discussion in  \cite{chawil01}), but there is far less discussion on the meaning of marginal value to a consumer. 

The meaning of \textit{value} is entirely clear in the case of TCLs and pool pumps:  this is reflected by the consumer who sets preferences on the thermostat or pool pump.   However,  there is no definition of \textit{marginal value} that is associated with power or any snap-shot definition of energy. 
   Let's use the insight from 
\Cref{s:consumers} to build a CE model that reflects the true utility of consumers for these types of loads.    To use the QoS model \eqref{e:SoC_ODE} we assume that there is an aggregator engaged with all of the consumers who ensures strict bounds on QoS.     It is simplest
to impose a cost $c_i\colon\Re\to\Re_+$  which is very large or infinite outside of the QoS interval.  
\rwm{does this imply we agree that it should be very large?
\\
sm2rwm:   if the aggregator has a contract with the customer to keep temperature within bounds, yes} 
The utility $\clU_D$ for the aggregator representing $M$ load classes is represented as the negative of the total cost:
\begin{equation}
\clU_D(P_1,\dots, P_M)  =  -   \int_0^{\clT}    \cX(x(t))   \, dt     =  -  \int_0^{\clT}  \sum_{i=1}^M   c_i(x_i(t)) \, dt
\label{e:clUD}
\end{equation}
In this formulation, the variable  $P_i $ now represents a \textit{function of time} on the interval $[0,\clT ]$ (e.g., $\clT=24$ hours in a day-ahead market).
For a given price process, the aggregator's optimization becomes
\begin{equation}
\begin{aligned}
\max \ \  & \clU_D(P_1,\dots, P_M)     -  \sum_i    \int_0^{\clT}   \varrho_t  d_i(t)  \, dt
   \\
 \text{subject to} \ \            &   \ddt x_i(t)   = -\alpha_i x_i(t)   - d_i(t),    \qquad 1\le i\le M
\end{aligned} 
\label{e:clUDopt}
\end{equation}
where  $d_i(t) = P_i(t) - P^b_i(t)$ where $P^b_i(t)$ is baseline power consumption from load class $i$ at time $t$  (such as the plot shown on the upper left hand side of  \Cref{f:WHTemp1} in the case of water heaters).   

Once again, this utility function is based on the average QoS of the population.  It is up to the aggregator to ensure that the QoS for each individual lies within pre-determined bounds.

   \begin{wrapfigure}[15]{L}{0.3\textwidth}

 \centering 
    \includegraphics[width=0.95\hsize]{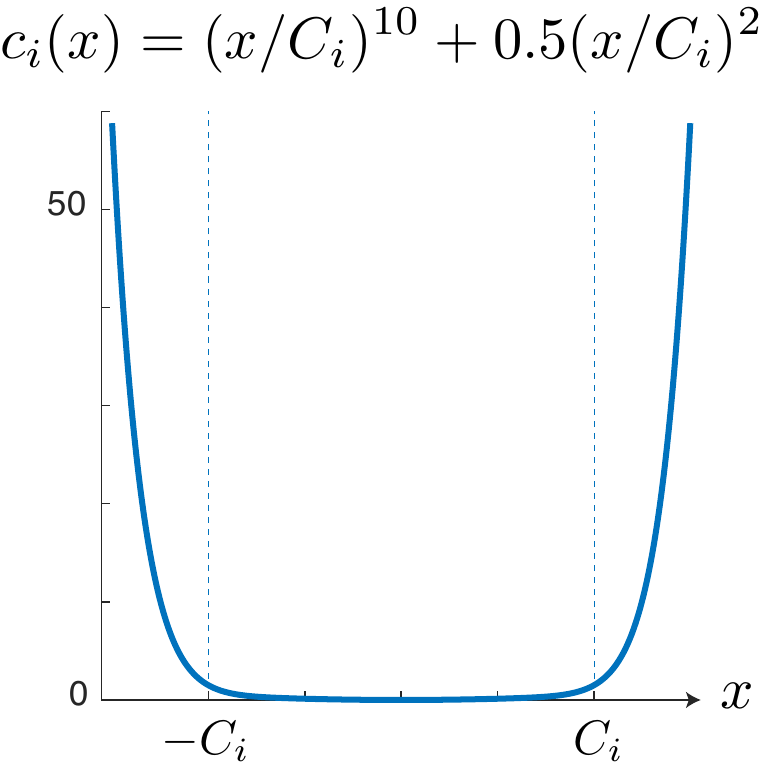}
	\caption{Cost of QoS violation.  }  	
 	\label{f:QoScost}
 
\end{wrapfigure}

\paragraph{Optimal CPP response}

The experiment described here is taken from the conference paper  \cite{matmoymeywar19} (see also \cite{JoelMathiasThesis21}).
Cost functions were chosen to be a high order polynomial, illustrated in   \Cref{f:QoScost}, 
with $\{ C_i  : 1\le i\le M \}$ capacity limits,   reflecting a ``snapshot'' QoS metric for the consumers. 
The choice of load quantities and characteristics in this experiment was based on a realistic population of loads in the CAISO area.

Consider a typical CPP event in which the price $\varrho_t$  takes on just two values, with a higher price occurring for only ninety minutes over the 24 hour period, and the increase is   by just 10\%.    The optimal solution to \eqref{e:clUDopt} shown in \Cref{f:boom} is obviously not the desired outcome from the point of view of the BA.    The 10~GW surge in power  followed by the instantaneous 13~GW drop  would cause blackout in any region of the U.S.

If not for the blackout, the consumers would be satisfied: air conditioning  and refrigerators turn off for between 15 and 30 minutes, and then return to the baseline power consumption.  Swimming pools and   water heaters respond approximately as predicted, but not smoothly:  regardless of the price increase,  these loads turn off for the entire 90 minute period since this results in no loss of QoS  (in terms of maintaining cleaning or temperature).

This is the biggest source of disaster in terms of load control:  all the loads turn off, regardless of the magnitude of price increase.    The most alarming outcome observed in this experiment is the ``pre-bound'' and not the rebound which is   commonly discussed in the literature  \cite{qdr2006benefits,Luetolf2018ReboundEO}.  This phenomenon is not the product of an eccentric model:  it should be obvious that  these   loads will get ready for a forecast CPP event by consuming extra power when prices are lower,  and equally obvious that they will turn off at the onset of a price surge.

In conclusion, the CPP pricing model fails for these deferrable loads.   

 \subsection{Dynamic Competitive Equilibrium Outcome}

We turn next to the question:  \textit{what does the CE solution look like based on a realistic utility function for consumers?}   To obtain an answer,  it is necessary to consider the utility function for generation, and then characterize   an equilibrium.

A dynamic CE model   involves $M+1$  players in the simplified model in which each of the $M$ load classes is represented by a single aggregator.  In addition, it is assumed that there is a single supplier (in practice, a single class of suppliers) that provide traditional generation $\bfmg$.   
It is assumed here that the utility function for   the supplier is the negative of cost, based on 
a convex function of generation and   its derivative, $\clU_S(g,\dot g) = -c_g(g)  -c_d(\dot g)$.
We also include inflexible load, denoted $\ell(t)$ at time $t$, so the supply=demand constraint becomes
\begin{equation}
g(t)  = \ell(t) +  P_\sigma (t)      
	=
	 \ell(t) +  P^b_\sigma (t) +  	d_\sigma (t)  
\label{e:supply=demand}
\end{equation}
where we have introduced the notional convention $d_\sigma = \sum_i d_i$ for any   $\{d_i\}$.

There is a price process  $\varrho^\star_t$ that forms a competitive equilibrium, defined as a solution to the following $M+1$ optimization problems
\begin{equation}
\begin{aligned}
 \bfmd_i^\star      & \in   \argmax_{d_i}     \int_0^\clT    \clU_{D_i}(x_i(t))    -  \varrho^\star_td_i(t) \, dt ,
\\
\bfmg^\star & \in \argmax_g   \int_0^\clT  \clU_S(g(t) ,\dot g(t))    +\varrho^\star_t  g(t)   \, dt  .
\end{aligned} 
\label{e:greedyPlayers}
\end{equation}
with  $ \clU_{D_i}(x_i(t))   =-   c_i(x_i(t)) $,   
and such that the constraints \eqref{e:SoC_ODE},  \eqref{e:supply=demand}, hold for the optimizers.
  
       To identify the equilibrium price, we follow the steps of \cite{wankownegshamey10} which is a dynamic extension of the standard textbook construction  \cite{maswhi1995}:  we posit an optimization problem whose objective is the sum of utilities,  and the price is then the Lagrange multiplier associated with the supply=demand constraint \eqref{e:supply=demand}.

The optimization problem is known in economics textbooks as
the \textit{social planner's problem} (SPP).   For the model considered here, the SPP can be expressed as  
\begin{gather}
\label{e:CEdyn}
\max_{g,d}  
 \int_{0}^{\clT}  \Bigl\{  
 		 \clU_S(g(t) ,\dot g(t))  +  \sum_{i=1}^M  \clU_{D_i}(d_i(t))   
		 		\Bigr\}  \, dt 
\end{gather} 
To simplify notation we transform to a minimization problem, and make the constraints explicit:
\begin{subequations}
\begin{align}
\textbf{SPP:}
\qquad \qquad
\min_{g,d}   \ \ &
	\int_0^{\clT}    \big[   c_g(g(t) ) + c_d(\dot g(t))  + \cX(x(t)) \bigr] \, dt 
 \label{qp19}
 \\
 \text{subject to} \ \ &
	P^b_\sigma (t) =g(t) -   d_\sigma (t)  - \ell(t)
 \label{e:balancecons}
	\\
	&
 \ddt {x}_i(t)  =- \alpha_i x_i(t) + d_i(t)
\label{e:soccons}%
\end{align}%
\label{e:CEconstraints}%
\end{subequations}
where  $\cX$ defined in \eqref{e:clUD},  and
  $x(0), d(0) \in \Re ^M $  are  given.    
    Let $\lambda_t$  denote the multiplier associated with  the constraint  \eqref{e:balancecons}, and consider the Lagrangian:
\[
\begin{aligned}
\clL(g,d,\lambda) = 
	\int_0^{\clT}  &  \big[   c_g(g(t) ) + c_d(\dot g(t))  + \cX (x(t)) \bigr] \, dt 
	\\
	&
	+	\int_0^{\clT}   \lambda_t  \big[    P^b_\sigma(t)  - g(t) +    d_\sigma (t)  + \ell(t) \bigr] \, dt  
\end{aligned} 
\]
The minimization of the Lagrangian is known as the dual functional:
\[
\varphi(\lambda)    \eqdef \min_{g,d} \clL(g,d,\lambda)   
\]
where the minimum is subject to the dynamics \eqref{e:soccons} and the given initial conditions.  The dual function admits a representation as a \textit{Lagrangian decomposition}:
\begin{equation}
\begin{aligned}
\varphi(\lambda)   &  =  \min_g  \int_0^{\clT}     \big[   c_g(g(t) ) + c_d(\dot g(t))  - \lambda_t g(t)     \bigr] \, dt 
            \\
	&\quad +\sum_i  \min_{d_i}  	\int_0^{\clT}   \bigl[   c_i (x_i(t))  +\lambda_t d_i(t)       \bigr] \, dt 
	\\
		&\qquad\qquad\qquad +   	\int_0^{\clT}   \lambda_t  [    P^b_\sigma(t)  + \ell(t)]  \, dt  
\end{aligned} 
\label{e:CEdecomposition}
\end{equation}
This decomposition is the main ingredient in the welfare theorems of  competitive equilibrium theory  (see \cite{maswhi1995} for similar decompositions in static economic models,  and the survey \cite{wannegkowshameysha11b} for dynamic CE models).   
The proof of \Cref{t:dualCE} follows from the definitions:   

\begin{proposition}
\label[proposition]{t:dualCE}
Suppose that there is an optimizer $\lambda^\star$ for the dual functional:  $ \varphi(\lambda^\star) = \max_\lambda \varphi(\lambda)$,   and respective optimizers  
\[
\begin{aligned}
d_i^\star	&\in  \argmin_{d_i}  	\int_0^{\clT}   \bigl[   c_i (x_i(t))  +\lambda^\star_t d_i(t)       \bigr] \, dt  \,,\qquad 1\le i\le M
\\
g^\star  &  \in  \argmin_g  \int_0^{\clT}     \big[   c_g(g(t) ) + c_d(\dot g(t))  - \lambda^\star_t g(t)     \bigr] \, dt 
\end{aligned} 
\]
Then \eqref{e:greedyPlayers}  is satisfied, so that   $\varrho^\star = \lambda^\star$ defines a price which results in a competitive equilibrium.
\end{proposition}


%

A similar optimization problem arises in  \cite{matmoymeywar19}, which concerns optimal control formulations for  resource allocation   (including a brief discussion on potential market implications).    Theorem~3.1 of \cite{matmoymeywar19} implies the following relationship between QoS and price:  for each $i$ and $t$,
\begin{equation}
-  c'_i\,   (x^\star_i(t))  =  \alpha_i \varrho^\star_t -  \ddt \varrho^\star_t 
\label{e:CEmarginalCost}
\end{equation}
The left hand might be interpreted as the \textit{marginal value} for players of class $i$.    
This conclusion is very distant from classical economic thinking, in which marginal value is equal to the equilibrium price.    This relationship between price and marginal value also implies that $\varrho^\star_t $ must be   continuous, which in particular rules out the price signal used in CPP today.

\begin{figure}[ht]
\Ebox{.65}{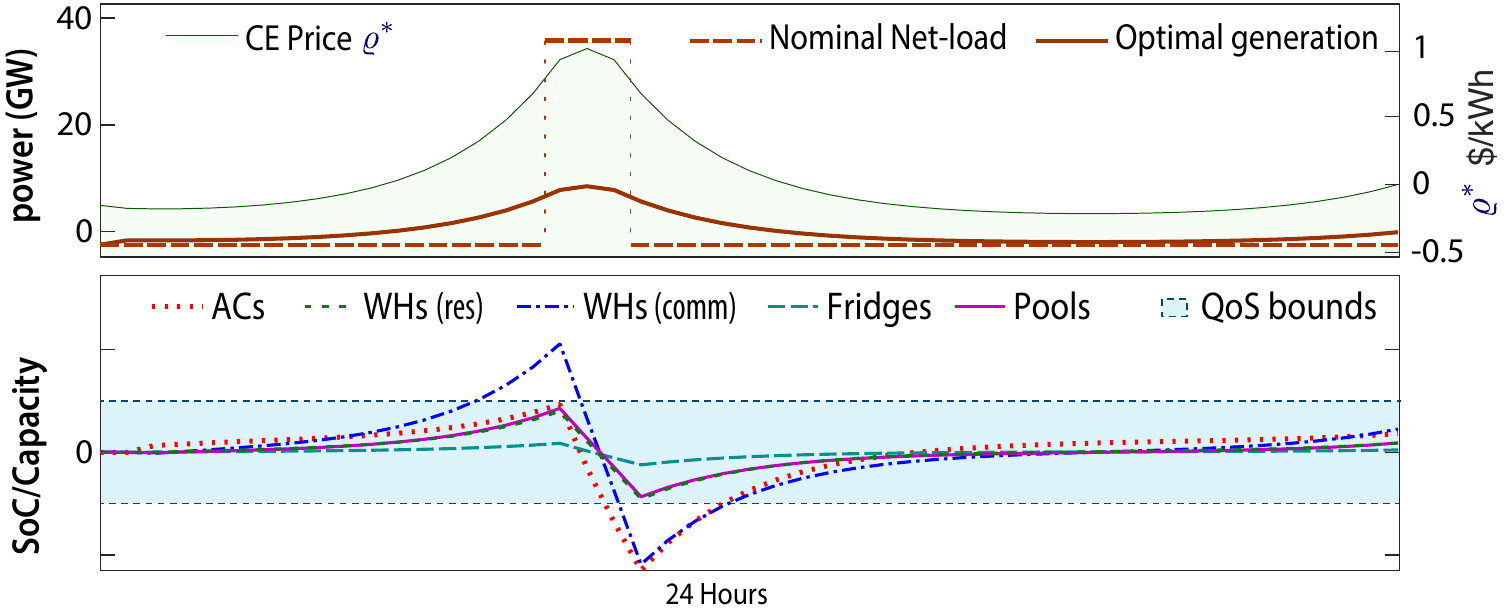}

\caption{Competitive Equilibrium Solution, with nominal net-load   $\ell(t) + P^b_\sigma(t)$ shown.}
\label{f:QuadSocCapacity_Policy} 
\end{figure}

Let's revisit the experiment leading to \Cref{f:boom}:  rather than impose a surge price,   we  consider a period of time in which netload is far greater than normal,  and compute the resulting CE price.       The emergency event occurs for 90 minutes as before, during which time the net-load increases by 40~GW. 
Based on the same cost/utility functionals used to create the data in \Cref{f:boom},  the CE price along with generation and consumer responses were computed.  The results are summarized in the plots shown in   \Cref{f:QuadSocCapacity_Policy}.      The price is very smooth, as anticipated from the identity \eqref{e:CEmarginalCost} and the fact that generation is subject to ramping costs.     
The loads ramp up and down gracefully, resulting in a highly smooth optimal trajectory for generation.


Observe that all the VES ``charge" prior to the increase in netload, so that they can ``discharge" during the increase. In other words, the ACs pre-cool and the WHs pre-heat in anticipation of the increase.   The commercial hot water heaters experience a loss of QoS for under one hour near the start of the scarcity event (the water is over-heated),   
and both commercial hot water heaters and air conditioning  sacrifice QoS for a similar period of time at the end of the event (the water is under-heated and the air conditioning is under-cooling).   This is a product of design. 
If  tighter bounds on QoS are required, then the cost function must be modified to reflect this,  or additional resources are required to maintain smooth output from thermal generation.  

%
  
%



\subsection{Summary \&\ Reality Check}

\Cref{t:dualCE} makes clear how dynamics can impact  markets:
if it is assumed that  the preferences of the consumer can be approximated by a concave function of power consumption, then we arrive at the dynamic pricing models that have been advocated for the past two decades.  We come to an entirely different conclusion upon recognizing that this model does not come close to reflecting user preferences for the majority of electric loads.  

It might be possible for a reader to see \Cref{f:QuadSocCapacity_Policy}, and declare that power pricing is solved by \Cref{t:dualCE}.   This would be a misreading of the proposition, which is meant to show that our conclusions are highly sensitive to modeling assumptions.  How do we know that aggregators and consumers  will behave exactly as our rationality assumptions predict?    Moreover,  a core ingredient of the theory of efficient markets is that prices must be  \textit{discovered} via some mechanism.   We feel it would be prohibitively difficult to design a real-time market that could discover the dynamic price described in the proposition.

Competitive equilibrium pricing models also make the fatal mistake of assuming that fixed costs are ``sunk'' and therefore not needed in the problem formulation. 
 This is a simplification imposed in \Cref{t:dualCE}  for tractability, and to compare our conclusions with all of the other papers on marginal cost pricing that impose the same simplification.  
 This is in fact far more than a simplification: it is a gross distortion of reality.  Evidence of this is that in these models, standard arguments show that average price is approximately average marginal cost   \cite{chomey10,wannegkowshameysha11b,negetal17}.   Obviously, generation companies are not simply optimizing over a 24-hour period based on cost of energy and ancillary services,  and it should also be obvious that their average marginal cost is only a fraction of the revenue required to maintain a healthy business \cite{moye2021}.

The proposition ignores the important impact of risk and uncertainty.   
Soon after the February 2021 crisis at ERCOT, 
a panel consisting of many architects of marginal cost pricing for power systems met to discuss the outcome 
 \footnote{\url{https://www.youtube.com/watch?v=Wz3172SIAfI}}.   All agreed that the actors (including consumers) should have hedged to avoid financial crisis.   This may be true, but the architects of these markets did not account for risk in their market designs.    Competitive equilibrium theory falls flat with the introduction of risk, as made clear in \cite{kaw05} and the references therein.

%
  
Moreover, our goals go far beyond efficiency.
 Given the enormous value of reliability,  we need planning and reliable control  mechanisms that do not depend on our assumptions about the rationality of all of the competing agents.  Reliability requires control loops, just as in telecommunications.   
 The cost of failure is far higher in the case of our power grids, in part because the recovery process is much more complex and time-consuming.

Finally, \textit{our energy grid is an interconnection of overlapping complex networks,  built on a myriad of dynamical systems}.  
It is difficult to envision marginal cost pricing as an effective control strategy in such a complex network.


\section{Planning and Operations in a System of Systems}

This section goes beyond the technical but common-sense conclusions of \Cref{s:idiots}, providing a critique of marginal cost pricing based on the authors' combined 60 years of   expertise in the power industry.     Many of our conclusions are based on a single observation:  designing and operating the power grid while ensuring reliability is a complex endeavor that requires rigorous long-term and short-term planning as well as advanced control techniques.

%
%
%



\paragraph{Marginal cost pricing fails to drive long-term investment}
Installing generators or power lines and supporting devices takes  years for the bulk system, and months for smaller assets.  For example, adding an operational transmission line (or a power plant) can take up to 10 years in California. The timeline is shorter in Texas, but generally not less than 3 years.  Deregulation was intended to drive a competitive market and simultaneously incentivize creation of more efficient assets.  This deregulated organization worked for a number of years when our grid had ample capacity to sustain reliability. 
 Over time, with energy as the commodity driving this market, and pricing as settlement, the system was driven towards thin margins of reserves and asset maintenance in many parts of the U.S., leaving the system skating on thin ice.
%
%


Deregulation has created a system of siloed stakeholders and largely siloed responsibilities.  Although these entities have complex processes and systems for handing off information to coordinate activities, we do not currently have the technology or transparency necessary to plan or operate  the grid as a holistic system. The most significant disconnect today is between transmission and distribution on both the grid and the resource level.  \textit{This is an inefficient means to operate an energy grid, lacking the coordinated metrics to maintain end to end reliability in planning and operation.}

In the case of the   recent 2021 power systems crisis at ERCOT, the following observations are almost self-evident (although the causes and implications of the crisis need to be studied in detail):  
\begin{romannum}
	

\item  	The market design was not adequate to incentivize the stakeholders to invest in system capacity, maintenance, winterization, or reliability. This lead to massive, widespread vulnerabilities across multiple layers of the grid: generation capacity, transmission networks, distribution networks, contingency reserves, etc.

\item Several stakeholders buckled under high electricity bills resulting in missed payments to grid operators and utility companies, thereby leading to several cases of bankruptcy. This is a clear indication that surge pricing did not mitigate the shortcoming of a market designed around the marginal cost of power.
\end{romannum}

\paragraph{Operating reserves are not planning reserves}
A more subtle but critical flaw in current pricing models is related to how balancing authorities prescribe \textit{operating reserves}, a flaw which was vividly exposed  during the recent ERCOT crisis.  During hours when service to thousands of customers had been curtailed (but operating reserves were sufficient), the market-wide LMP was as low as \$1,200/MWh. Recognizing the obvious disconnect between the VOLL-based scarcity pricing construct and actual market conditions, the Public Utility Commission of Texas took the “unprecedented step" of requiring ERCOT to set market prices at \$9,000/MWh. They concluded that “if customer load is being shed" the price should be \$9,000/MWh. ERCOT complied with this directive and kept prices at that level until it declared an end to the power emergency on Friday, February 19, 2021  \cite{wsj2021a}.    Those loads that are shed are providing ancillary service to the grid, ensuring sufficient operating reserves.  It seems obvious that these customers should be receiving  compensation during this crisis, and not the generation companies.  The generators should be receiving compensation in advance for \textit{planning reserves} so as to minimize the frequency of costly load shedding events.
\jjm2021{Mar 10: not sure this is needed,
 but that is not how the market is working.  Although the data is not available yet, marginal cost of power and surge pricing did not successfully incentivize generators to utilize available market products to hedge market risks}

\paragraph{We can reliably control complex systems}   It would appear that all of these issues are exacerbated by an increasingly complex, uncertain grid which involves many activities, including balancing resources, forecasting,
 load switching, integration of renewables, monitoring transmission, 
	managing fuel supply infrastructure,
scheduling planned generator downtime, and managing forced outages.
This is a lot to oversee, even without  accounting for the unpredictable behavior of humans interacting with these  complex systems.  
Theory and recent events make clear that a simplistic market design will likely exacerbate these complexities.   On the other hand, our success in managing complex supply chains, transportation systems and communication networks make clear that reliability can be greatly improved by taking advantage of a rich toolkit from decision and control theory.

\jjm2021{We need either citations or need to make this more precise: "real-time markets make short-term planning to balance the resources on the grid, from the generators, balancing authority through the transmission/distribution all the way to the consumer prone to delays, inefficiencies, lack of synchronization, as well as deficiencies in systems, processes and responsibilities."}


\section{Alternatives}
\label{s:alt}

We turn next to policy recommendations and questions that must be answered as we improve our energy infrastructure.

 \subsection{Reliability System Operator as a Central Planner}   
 \label{s:RSO}
 
 Every successful industry performs long-term planning in the face of risk.  The CEO of Delta Airlines works with a team of experts to plan the next fleet of airplanes.  Markets play an important role when airplane manufacturers bid for contracts,  and the bargaining between suppliers and consumers leads to substantial innovation.  \textit{Why do we fear the CEO model in power systems operations?}

Under the recommended solution, the reliability system operator (RSO) acts  as a central planner and develops an optimal resource expansion plan across the entire market footprint using traditional planning techniques and methodologies. However, the proposed process incorporates a level of sophistication and complexity similar to that developed for the operation of organized markets.  
The resource providers will be selected on the basis of a competitive techno-economic process, and the prices for these resources will be reflected in the long-term contracts signed with the successful suppliers. The prices would incorporate both sunk and operational costs.\jjm2021{See edits to the above paragraph. Also added last sentence.}

Similar to approaches still utilized in the traditionally regulated jurisdictions in the Western and Southeastern U.S., the RSO will conduct system-wide planning in order to fulfill the following responsibilities, many of which are currently performed by Regional Transmission Operators (RTOs). \jjm2021{RTO replaced by RSO everywhere below, and deleted some bullet points to avoid repetition}
\begin{itemize}
	\item Reliability requirements are established and tracked. 
  	\item Future load requirements are forecast.
  	\item Planned generation, transmission, and potential fuel supply assets are identified and incorporated.
  	\item Long-term analyses are performed by the RSO that identify:
    \begin{itemize}
	  	\item the amount of capacity needed,
	  	\item the desired location of the capacity considering existing and potential transmission, and natural gas (or other fuel supply) infrastructure, and
	  	\item the additional technology needed to provide the capacity. 
  	\end{itemize}
\end{itemize}

The proposed approach would center around these traditional operations, but evolve to introduce avenues for competitive markets that emphasize reliability services and contracts evaluated under techno-economic performance indices. Under the proposed approach, following the development and agreement on a plan for the market, competitive auctions are held to select suppliers to build, own, and operate the resources. The RSO signs contracts with the successful bidders for the purchase of the capacity and associated energy under long-term (e.g., 10+ years) agreements, with such contracts having strict performance and financial guarantees (techno-economic performance indices).

Similarly, load-serving entities (LSEs) would enter into long-term contracts with the RSO for at least a portion of their requirements.  For the remaining portion of the load serviced by the LSE, it is compensated on the basis of its load-share ratio. This approach will also facilitate investments in demand-side solutions, which are discussed next.

\subsection{Demand Dispatch}
\label{s:sanity}
 
The DoE technical report \cite{qdr2006benefits}  was written in response to the Energy Policy Act of 2005, which asserted that it is the policy of the United States to encourage \textit{time-based pricing and other forms of demand response}.  The report presented a formal definition of demand response that is now widely accepted:
\begin{quote}
[Demand response refers to] ``changes in electric usage by end-use customers from their normal
consumption patterns in response to changes in the price of electricity over time, or to incentive
payments designed to induce lower electricity use at times of high wholesale market prices or
when system reliability is jeopardized''  \cite{qdr2006benefits}.
\end{quote} 
Recall from \Cref{s:history}   that demand response programs can be classified into two categories,  distinguished by the methods in which pricing is implemented:  the direct load control implemented by FPL is of the first category, in which the consumer receives an annual reward to allow interruption of service.  The second category is truly real-time control through prices that vary perhaps hour by hour, now known as the \textit{prices-to-devices} approach to load control. 
References  \cite{qdr2006benefits,jor19}   contains   history on both approaches.

\paragraph{Alternatives}
It should be clear from the theory and examples in this paper that ``prices-to-devices'' is not a workable approach to load control for the majority of electric loads.   

The 2010 article \cite{brolureispiwei10} introduced the term  \textit{demand dispatch} to describe a  third approach to consumer engagement.  
While there are some similarities with direct load control, it is anticipated that demand dispatch will be based on far more distributed intelligence.   
Returning to water heaters, it is the thermometer at the load that is measuring temperature, which is what determines QoS.  It is absolutely essential then to create intelligence at that load so that   bounds on QoS are maintained, while simultaneously providing service to the grid.

The article makes a much larger observation: much of the power economics literature is concerned with reducing consumption peaks through price signals, which is the thinking behind critical peak pricing (that this paper is warning against). Surely, through intelligent control design, we can obtain an enormous range of grid services, beyond just peak reduction, by enlisting an army of water heaters, refrigerators, commercial and residential HVAC, etc.

The question of how to organize an ensemble of loads to obtain reliable VES has been a focus of research since the late 1970s \cite{sch78,schFAPER80,malcho85},  and saw a resurgence of interest in the past decade:  see the surveys included in  \cite{IMA18}, in particular \cite{cheche17b,almesphinfropauami18,chehasmatbusmey18,moymey17a}.   Notable articles and theses include   \cite{johThesis12,haomidbarey12,linbarmeymid15,YueChenThesis16, matbusmey17, cammatkiebusmey18,bencolmal19}.
It seems now that distributed control for reliable VES is a mature discipline that is ripe for application.  

\textit{Next Steps}:
\begin{itemize}
\item
 Field testing is required to discover potential gaps between theory and application.   
 
 \item
 Research is required on the benefits and risks at the distribution layer.  For example,   VES is surely useful for voltage regulation.   Can loads provide voltage support and grid level services simultaneously?    

 \item
There is the question of VES availability during extreme events.    
During the 2011 heat wave in Texas, many AC units were operating at full power during much of the day.   
The capacity from VES will be reduced in such cases, as well as the capacity from traditional generation:  how can a generator provide ramping services when it is operating at maximal output?  
This means we still require standby resources for the most extreme  scenarios (in this case, once in a decade),   as well as more sophisticated approaches to controlled blackouts.

%
\end{itemize}
We recommend a guiding principle in any demand dispatch design:  \textit{do not implement if you cannot simulate}.    That is, if a control design or market mechanism is too complex to be properly evaluated, then it is probably too risky to install on our power grid.  


\section{Conclusions}
\label{s:con}

The control community knows that a complex engineering problem is not solved simply because a theorem has been proven.   Testing is required through simulation, and then in the field.  When the control system fails,   assumptions are revisited and the design is improved.     
Mathematical economics often avoids these critical next steps.  We are told that theory predicts price discovery, but 
we can never predict the genius of the market.  Consequently, no simulation can ever predict what magic will unfold.

However, it should be obvious that this \textit{genius} is a two-edged sword.   In the case of critical infrastructure, we do not want the uncertainty that comes with self-interested agents that seek to maximize profits,  without motivation to maintain reliability.   
\jjm2021{Mar 10 Looks la bit hubristic: As industry leaders and researchers, }We appreciate all the innovations that come with the marketplace, but we know from examples such as Apple and Delta Airlines that the term \textit{market} is not synonymous with \textit{marginal cost pricing}.

It is hoped that this paper will accelerate the evolution away from short-term pricing models, 
and towards a greater appreciation of long-term planning and business models based on techno-economic metrics.   We will not lose the genius of the market.  The innovators will be encouraged to think years ahead into the future instead of just the next 24 hours.

\bibliographystyle{abbrv}

\def\cprime{$'$}\def\cprime{$'$}


\clearpage

 \section{Biographies}
 
 \paragraph{Hala Ballouz} is a lead consultant in the electric power industry since 1991.  She is the Owner and President of Electric Power Engineers, Inc.  and a registered Professional Engineer (P.E.) with expertise focused on power systems and energy market studies.  She is also the founder of Utility Plexus Engineering, a firm focused on innovation on software + consulting for modernizing utility grids.  Ms. Ballouz has been an advisor nationally and internationally on the integration of renewables and smart grid development.  She is an innovator in the holistic study of planning and operation of transmission and distribution (T\&D), and integration of generation, load, and distributed resources, including implementation of DERMS and Advanced Power Management Platforms. She   founded GridNEXT, while President of the Texas Renewables Energy Association (TREIA) and established the first of a kind series of conferences focusing on grid modernization to further the integration of renewables. She is currently leading the integration of smart grids into T\&D planning and operation through innovative software development and application; with her vision, she set a founding cornerstone at Texas A\&M University's Smart Grid Center to integrate ENER-i distribution system analytics \& DER platform to fuel innovation in holistic electrical grid planning.  She is a recognized leader in power systems as well as in empowering women in the energy industry. Through an excellent team of power systems engineers, her company is on a path of growth to exceed 40\% annually. Hala Ballouz  holds a Bachelor of Science and a Master of Science degree in Electrical Engineering.  
She was the recipient of the 2020 \textit{Woman in Power} Award.

 \paragraph{Joel Mathias} is currently pursuing a Ph.D.~degree in electrical and computer engineering at the University of Florida. He received the M.S.~degree in Electrical and Computer Engineering at the University of Florida, Gainesville, FL, USA, and the B.Eng.~degree in Electronics and Communications Engineering from the University of Mumbai. He was a research associate with the Tata Institute of Fundamental Research, Mumbai, India, from 2011--2012, and his industry experience includes Electric Power Engineers, Inc., Austin, TX, USA, and Tata Consultancy Services, Mumbai, India. For his graduate studies in the USA, he has received the  J.N.~Tata Endowment for the higher education of Indians, the Lady Navajbai Tata scholarship, and achievement awards from the University of Florida. His research interests are in stochastic processes, control theory, and optimization, with applications in the area of smart power grids.

 \paragraph{Sean Meyn}  is the Robert C. Pittman Eminent Scholar Chair at the Department of Electrical and Computer Engineering at the University of Florida, director of the Laboratory for Cognition and Control, and Inria International Chair at Inria, France.  He is a Fellow of IEEE.  He joined Florida in 2012, following 20 years as a professor of ECE at the University of Illinois. Following his BA in mathematics at UCLA, he moved on to pursue a PhD with Peter Caines at McGill University. During 2018 he spent his sabbatical at the Simons Institute at Berkeley in the spring,  and Inria Paris, and Ecole Polytechnique during the fall. His interests span many aspects of stochastic control, stochastic processes, information theory, and optimization. For the past decade, his applied research has focused on engineering, markets, and policy in energy systems.  His monograph “Markov Chains and Stochastic Stability” (co-authored with R.L. Tweedie) is a standard text in algorithm analysis and many other fields, and was awarded the 1994 ORSA–TIMS best publication award.  In 2015 he and Prof.~Ana Busic received a Google Research Award, “to foster collaboration on demand dispatch for renewable energy integration”.  He is a CSS distinguished lecturer on topics related to both reinforcement learning and energy systems, and is active in organizing programs on these topics.  Most recently, he was co-organizer of the Simons Institute program on the Theory of Reinforcement Learning, Aug. 19--Dec. 18, 2020.

\paragraph{Bob Moye} joined Tyr Energy in December 2019 to establish and grow the company’s energy management services business, focusing on both Tyr-owed as well as third-party assets.  Prior to Tyr, Mr. Moye worked for 14 years developing and managing Rainbow Energy Marketing Corporation’s energy management group.  Over that time, arrangements with 27 customers were executed representing 34 separate assets and 7 electric utility systems consisting of generation and load, and capacity under management totaling nearly 8,000 MW.    Prior to joining Rainbow, Mr. Moye worked for The Energy Authority and was responsible for the growth of this newly created municipal energy management company.   During his tenure, he also managed all long-term marketing and risk management, and all short- and long-term physical and financial trading. Mr. Moye has Bachelor of Science in Electrical Engineering, Master of Engineering, and Master of Business Administration degrees.  He is currently a candidate for a PhD in Electrical Engineering at the University of Florida.  He will defend his dissertation in spring 2021, entitled \textit{Resource Investments in Organized Markets:  A Case for Central Planning}.

\paragraph{Joseph Warrington} is a Senior Staff Research Engineer at Home Experience LLC, Cambridge UK. Until December 2019 he was a Lecturer and Senior Scientist in the Automatic Control Lab (IfA) at ETH Zurich, Switzerland. His Ph.D.~is from ETH Zurich (2013), and his B.A.~and M.Eng.~degrees in Mechanical Engineering are from the University of Cambridge (2008). From 2014-2016 he worked as an energy consultant at Baringa Partners LLP, London, UK, and he has also worked as a control systems engineer at Wind Technologies Ltd., Cambridge, UK, and privately as an operations research consultant. He is the recipient of the 2015 ABB Research Prize for an outstanding PhD thesis in automation and control, and a residential Simons-Berkeley Fellowship for the period January-May 2018. His research interests include dynamic programming, large-scale optimization, and predictive control, with applications including power systems and transportation networks.

\end{document}